\documentclass{article}
\usepackage{amsmath,amssymb,mathtools}
\usepackage{tikz,xcolor}
\usetikzlibrary{backgrounds}
\usepackage{graphicx}
\usepackage{subfigure}
\usepackage{multirow}
\usepackage{tikz}
\usetikzlibrary{arrows,shapes,trees,calc}

\def\RR{\mathbb R}

\def\b1{\mathbf 1}

\newcommand{\raf}[1]{(\ref{#1})}

\newcommand{\hide}[1]{}

\newtheorem{remark}{Remark}

\renewcommand{\beta}{{\cal B}}

\title{A four-person chess-like game
without Nash equilibria in pure stationary strategies}

\begin{document}

\maketitle

\author{
\noindent
Vladimir Gurvich,
% \newline
\thanks{RUTCOR and RBS, Rutgers University, 100 Rockafeller Road, Piscataway, NJ 08854;
\newline
gurvich@rutcor.rutgers.edu; vladimir.gurvich@gmail.com}
\and
\newline
\noindent
Vladimir Oudalov,
% \newline
\thanks{80 Rockwood ave \#A206, St. Catharines ON L2P 3P2 Canada;
% \newline
oudalov@gmail.com}

\abstract{
In this short note we give the example of
a four-person finite positional game with perfect information
that has no positions of chance and
no Nash equilibria in pure stationary strategies.
The corresponding directed graph has
only one directed cycle and only five terminal positions.
It remains open:
\begin{itemize}
\item{(i)}
if the number  $n$  of the players can be reduced from $4$ to  $3$,
\item{(ii)}
if the number $p$ of the terminals can be reduced from  $5$  to $4$, and most important,
\item{(iii)}
whether it is possible to get a similar example in which
the outcome  $c$ corresponding to all
(possibly, more than one) directed cycles is
worse than every terminal for each player.
\end{itemize}
Yet, it is known that
\begin{itemize}
\item{(j)}
$n$  cannot be reduced to  $2$,
\item{(jj)}
$p$  cannot be reduced to  $3$, and
\item{(jjj)}
there can be no similar example in which
each player makes a decision in a unique position.
\end{itemize}

\smallskip

}

\section{Introduction}
\label{s0}

In this paper we answer in the negative to
one of the two open questions suggested in \cite{BEGM11}.
This negative result was conjectured in \cite{GO14}.
We refer the reader to these two papers for definitions of
the chess-like and backgammon-like games, and of of Nash-solvability,
which we recall here just briefly.

The {\em backgammon-like} and {\em chess-like} games are finite
positional $n$-person games with perfect information, which
can and, respectively, cannot have random moves.

More precisely, such a game is modeled by a finite directed graph
(digraph)  $G = (V,E)$, whose vertices are partitioned into
$n+2$  subsets: $V = V_1 \cup \ldots \cup V_n \cup V_T \cup V_R$.
A vertex  $v \in V_i$  is interpreted as a {\em position}
controlled by the player  $i \in I = \{1, \ldots, n\}$, while
$v \in V_R$  is a positions of chance,
with a given probabilistic distribution on the outgoing edges;
finally, $v \in V_T = \{a_1, \ldots, a_p\}$  is
a {\em terminal position} (or just a {\em terminal}, for short), which
has no outgoing edges.
Furthermore, a directed edge  $(v, v')$  is interpreted as
a {\em move} from the position  $v$  to  $v'$.
We also fix an initial position $v_0 \in V \setminus V_T$.

A game is called {\em chess-like} if
it has no positions of chance, $V_R = \emptyset$.

The digraph  $G$  %%%of a backgammon- or chess-like game
may have {\em directed cycles} (dicycles).
Recall that positions may  be repeated in a backgammon or in a chess play.
We assume that all dicycles  of  $G$  form a unique outcome  $c$  of the game.
Thus, the set of outcomes is  $A = \{a_1, \ldots, a_p; c\}$.

\begin{remark}
In  \cite{BGMS07} a different approach was suggested
(for  $n=2$): each dicycle was treated as a separate outcome.
Anyway, our main example contains only one dicycle.
\end{remark}

To each player  $i \in I$  and outcome  $a \in A$  we assign a payoff
(called in the literature also a reward, utility, or profit)  $u(i, a)$
of the player  $i \in I$ in case the outcome  $a \in A$  is realized.
The corresponding mapping
$u : I \times A \rightarrow \RR$  is called the payoff
(reward, utility, or profit) function.

Since our main result is negative and
related to chess-like games, we could restrict ourselves and
the players to their strict preferences, instead of the real-valued payoffs.
A complete order  $o_i$  over  $A$
is called the {\em preference} of the player  $i \in I$; and let
$o = (o_1, \ldots, o_n)$  denote a {\em preference profile}.

A {\em backgammon-like game in the positional form}
is the quadruple  $(G, D, o, v_0)$, where
$G = (V,E)$  is a digraph,
$D: V = V_1 \cup \ldots \cup V_n \cup V_T \cup V_R$
is a partition of the positions,
$o = (o_1, \ldots, o_n)$  is a {\em preference profile}, and
$v_0$ is a fixed initial position.
% The game is called {\em chess-like} if
% there are no positions of chance, $V_R = \emptyset$.
The triplet  $(G, D, v_0)$  is called a {\em positional game form}.

\smallskip

To define the {\em normal form}
(of a chess-like game)
let us introduce the strategies.
A {\em (pure and stationary) strategy} of a player  $i \in I$ is
a mapping that assigns a move  $(v, v')$  to each position  $v \in V_i$.
(In this paper we restrict ourselves and the players
to their pure stationary strategies, so
mixed and history dependent strategies will not even be introduced.)
A set of  $n$ strategies  $s = \{s^i, i \in I\}$  is called
a {\em strategy profile} or a {\em situation}.
Each situation uniquely defines a {\em play}  $P(s)$
that begins in  $v_0$  and
either ends in a terminal  $a \in V_T$  or cycles.
In the latter case  $P(s)$  looks like a ``lasso":
it consists of an initial part and a dicycle repeated infinitely.
This is so, because a (pure stationary) strategy
assigns the same move whenever a position is repeated
and, hence, each situation  $s$  uniquely defines a move  $(v, v')$
in each non-terminal position  $v \in V \setminus V_T$.
Thus, we obtain a {\em game form}, that
is, a mapping  $g : S \rightarrow A$, where
$S = S_1 \times  \ldots \times S_n$
is the direct product of the sets
$S_i = \{s^i_1, \ldots, s^i_{k_i}\}$  of strategies of all players $i \in I$.
The {\em normal form} of a chess-like game
$(G, D, o, v_0)$  is defined as the pair $(g, o)$.

For the backgammon-like games each strategy profile  $s$  uniquely
determines a Markovian chain, which assigns to each outcome  $a \in A$
(a terminal or an infinite play) a well defined limit probability $p(s, a)$.
The payoff  $u(i, s)$  of a player  $i \in I$  in this situation  $s$
is defined as the expectation of the corresponding payoffs
$u(i, s) \sum_{a \in A} p(s, a) u(i, a)$.

\smallskip

A situation  $s \in S$  is called a {\em Nash equilibrium}
(NE) if for each player  $i \in I$  and
for each situation  $s'$  that may differ from  $s$
only in the coordinate  $i$, the inequalities 
$o_i(g(s)) \geq o_i(g(s'))$  and  $u(i,s)) \geq u(i, s')$ 
hold in case  of the chess- and backgammon-like games, respectively;
or in other words, if
no player  $i \in I$  can profit replacing
his/her strategy  $s^i$  in  $s$  by a new strategy  $s'^i$, provided
the  $n-1$  remaining players keep their strategies unchanged.
Let us remark that 
the equality  $o_i(g(s)) = o_i(g(s'))$ 
may hold only when  $g(s) = g(s')$, since  
the preference  $o_i$  is strict. 

\begin{figure}
\label{f1}
\begin{center}
\begin{tikzpicture}[>=stealth]
  \node at (canvas polar cs:radius=2.3in,angle=-90) [shape=circle,minimum size=1cm,draw=black,thick] (U4) {$u_4$};
  \node at (canvas polar cs:radius=2.3in,angle=-162) [shape=circle,minimum size=1cm,draw=black,thick] (W2) {$w_2$};
  \node at (canvas polar cs:radius=2.3in,angle=-18) [shape=circle,minimum size=1cm,draw=black,thick] (V2) {$v_2$};
  \node at (canvas polar cs:radius=2.3in,angle=54) [shape=circle,minimum size=1cm,draw=black,thick] (V3) {$v_3$};
  \node at (canvas polar cs:radius=2.3in,angle=126) [shape=circle,minimum size=1cm,draw=black,thick] (U3) {$u_3$};
  \node at (canvas polar cs:radius=3in,angle=-90) [shape=circle,minimum size=1cm,draw=black,thick] (A4) {$a_4$};
  \node at (canvas polar cs:radius=3in,angle=-162) [shape=circle,minimum size=1cm,draw=black,thick] (A5) {$a_5$};
  \node at (canvas polar cs:radius=3in,angle=-18) [shape=circle,minimum size=1cm,draw=black,thick] (A3) {$a_3$};
  \node at (canvas polar cs:radius=3in,angle=54) [shape=circle,minimum size=1cm,draw=black,thick] (A2) {$a_2$};
  \node at (canvas polar cs:radius=3in,angle=126) [shape=circle,minimum size=1cm,draw=black,thick] (A1) {$a_1$};
  \node at (canvas polar cs:radius=1in,angle=126) [shape=rectangle,minimum size=1cm,draw=black,thick] (U1) {$u_1$};
  \node at (canvas polar cs:radius=1in,angle=54) [shape=circle,minimum size=1cm,draw=black,thick] (U2) {$u_2$};
  \node at (canvas polar cs:radius=1in,angle=-90) [shape=circle,minimum size=1cm,draw=black,thick] (V1) {$v_1$};

  \draw [->] (U1) to (U2);
  \draw [->] (U1) to (U3);
  \draw [->] (U2) to (V1);
  \draw [->] (U2) to (V3);
  \draw [->] (V1) to (V2);
  \draw [->] (V1) to (U4);
  \draw [->] (V1) to (W2);
  \draw [->] (U3) to (V3);
  \draw [->] (U3) to (A1);
  \draw [->] (V3) to (V2);
  \draw [->] (V3) to (A2);
  \draw [->] (V2) to (U4);
  \draw [->] (V2) to (A3);
  \draw [->] (U4) to (W2);
  \draw [->] (U4) to (A4);
  \draw [->] (W2) to (U3);
  \draw [->] (W2) to (A5);

\end{tikzpicture}
\end{center}

\caption{
% \label{f1}
This figure represents our main example in the positional form.
\newline
Four players $I = \{1,2,3,4\}$  make decisions
in eight non-terminal positions
$u_1, v_1; u_2, v_2, w_2; u_3, v_3$, and $u_4$, respectively.
The subscript  is the number of the player
who controls the corresponding position.
\newline
The initial position is  $u_1$. 
There are five terminal positions  $a_j : j \in J = \{1,2,3,4,5\}$.
\newline
There is a unique dicycle  $c$ and, 
thus, the set of outcomes is  $A = \{a_1, a_2, a_3, a_4, a_5; c\}$.
\newline
The game has no NE  whenever the preferences  $o_i$  of the players $i \in I$
over the set of outcomes  $A$  agree with the following partial orders:
\newline
$O_1: a_2 > a_4 > a_3 > a_1 > a_5; \;\;$
%\newline
$O_2 : \min(a_1, c) > a_3 > \max(a_4, a_5) > \min(a_4, a_5) > a_2; \;\;$
\newline
$O_3 : \min(a_5, c) > a_1 > a_2 > \max(a_3, a_4); \;\;$
%\newline
$O_4 : \min(a_1, a_2, a_3, a_5) > a_4 > c$.
}
\end{figure}

\begin{table}
\label{t1}
\begin{tabular}{cc|cccccccc|cccccccc|}
\cline{3-18}
& & \multicolumn{8}{ c| }{$s^{4}_{1}$} & \multicolumn{8}{ c| }{$s^{4}_{2}$}
\\ \cline{3-18}
& & $s^{2}_{1}$ & $s^{2}_{2}$ & $s^{2}_{3}$ & $s^{2}_{4}$ & $s^{2}_{5}$ & $s^{2}_{6}$ & $s^{2}_{7}$ & $s^{2}_{8}$
& $s^{2}_{1}$ & $s^{2}_{2}$ & $s^{2}_{3}$ & $s^{2}_{4}$ & $s^{2}_{5}$ & $s^{2}_{6}$ & $s^{2}_{7}$ & $s^{2}_{8}$
\\ \cline{1-18}
\multicolumn{1}{ |c  }{\multirow{6}{*}{$s^{3}_{1}$} } &
\multicolumn{1}{ |c| }{$s^{1}_{1}$} & $c^{4}$ & $c^{4}$ & $a_{3}^{2}$ & $a_{3}^{23}$ &
 $a_{5}^{2}$ & $a_{5}^{2}$ & $a_{3}^{2}$ & $a_{3}^{23}$ &
 $a_{4}^{2}$ & $a_{4}^{23}$ & $a_{3}^{1}$ & $a_{3}^{3}$ &
$a_{4}^{24}$ & $a_{4}^{234}$ & $a_{3}^{1}$ & $a_{3}^{3}$
\\
\multicolumn{1}{ |c  }{}                        &
\multicolumn{1}{ |c| }{$s^{1}_{2}$} & $c^{4}$ & $c^{4}$ & $a_{3}^{23}$ & $a_{3}^{23}$ &
 $a_{5}^{2}$ & $a_{5}^{2}$ & $a_{3}^{23}$ & $a_{3}^{23}$ &
 $a_{4}^{23}$ & $a_{4}^{23}$ & $a_{3}^{13}$ & $a_{3}^{3}$ &
$a_{4}^{234}$ & $a_{4}^{234}$ & $a_{3}^{13}$ & $a_{3}^{3}$
\\
\multicolumn{1}{ |c  }{}                        &
\multicolumn{1}{ |c| }{$s^{1}_{3}$} & $c^{4}$ & $c^{4}$ & $a_{3}^{23}$ & $a_{3}^{23}$ &
 $a_{5}^{2}$ & $a_{5}^{2}$ & $a_{5}^{12}$ & $a_{3}^{23}$ &
 $a_{4}^{2}$ & $a_{4}^{23}$ & $a_{4}^{24}$ & $a_{3}^{3}$ &
$a_{4}^{24}$ & $a_{4}^{234}$ & $a_{4}^{24}$ & $a_{3}^{3}$
\\
\multicolumn{1}{ |c  }{}                        &
\multicolumn{1}{ |c| }{$s^{1}_{4}$} & $c^{4}$ & $c^{4}$ & $a_{3}^{23}$ & $a_{3}^{23}$ &
 $a_{5}^{2}$ & $a_{5}^{2}$ & $a_{3}^{23}$ & $a_{3}^{23}$ &
 $a_{4}^{23}$ & $a_{4}^{23}$ & $a_{3}^{13}$ & $a_{3}^{3}$ &
$a_{4}^{234}$ & $a_{4}^{234}$ & $a_{3}^{13}$ & $a_{3}^{3}$
\\
\multicolumn{1}{ |c  }{}                        &
\multicolumn{1}{ |c| }{$s^{1}_{5}$} & $c^{4}$ & $c^{4}$ & $a_{3}^{23}$ & $a_{3}^{23}$ &
 $a_{5}^{2}$ & $a_{5}^{2}$ & $a_{5}^{12}$ & $a_{3}^{23}$ &
 $a_{4}^{23}$ & $a_{4}^{23}$ & $a_{3}^{13}$ & $a_{3}^{3}$ &
$a_{5}^{12}$ & $a_{4}^{234}$ & $a_{5}^{12}$ & $a_{3}^{3}$
\\
\multicolumn{1}{ |c  }{}                        &
\multicolumn{1}{ |c| }{$s^{1}_{6}$} & $c^{4}$ & $c^{4}$ & $a_{3}^{23}$ & $a_{3}^{23}$ &
 $a_{5}^{2}$ & $a_{5}^{2}$ & $a_{3}^{23}$ & $a_{3}^{23}$ &
 $a_{4}^{23}$ & $a_{4}^{23}$ & $a_{3}^{13}$ & $a_{3}^{3}$ &
$a_{4}^{234}$ & $a_{4}^{234}$ & $a_{3}^{13}$ & $a_{3}^{3}$

\\ \cline{1-18}
\multicolumn{1}{ |c  }{\multirow{6}{*}{$s^{3}_{2}$} } &
\multicolumn{1}{ |c| }{$s^{1}_{1}$} &  $a_{1}^{3}$ & $a_{1}^{3}$  & $a_{3}^{23}$ & $a_{3}^{23}$ &
 $a_{5}^{12}$ & $a_{5}^{12}$ & $a_{3}^{2}$ & $a_{3}^{23}$ &
 $a_{4}^{24}$ & $a_{4}^{234}$ & $a_{3}^{1}$ & $a_{3}^{3}$ &
$a_{4}^{24}$ & $a_{4}^{234}$ & $a_{3}^{1}$ & $a_{3}^{3}$
\\
\multicolumn{1}{ |c  }{}                        &
\multicolumn{1}{ |c| }{$s^{1}_{2}$} &  $a_{1}^{3}$ & $a_{1}^{3}$  & $a_{1}^{1}$ & $a_{1}^{1}$ &
 $a_{1}^{3}$ & $a_{1}^{3}$ & $a_{1}^{1}$ & $a_{1}^{1}$ &
 $a_{1}^{1}$ & $a_{1}^{1}$ & $a_{1}^{1}$ & $a_{1}^{1}$ &
$a_{1}^{1}$ & $a_{1}^{1}$ & $a_{1}^{1}$ & $a_{1}^{1}$
\\
\multicolumn{1}{ |c  }{}                        &
\multicolumn{1}{ |c| }{$s^{1}_{3}$} &  $a_{1}^{3}$ & $a_{1}^{3}$  & $a_{1}^{1}$ & $a_{3}^{23}$ &
 $a_{5}^{12}$ & $a_{5}^{12}$ & $a_{5}^{12}$ & $a_{3}^{23}$ &
 $a_{4}^{24}$ & $a_{4}^{234}$ & $a_{4}^{24}$ & $a_{3}^{3}$ &
$a_{4}^{24}$ & $a_{4}^{234}$ & $a_{4}^{24}$ & $a_{3}^{3}$
\\
\multicolumn{1}{ |c  }{}                        &
\multicolumn{1}{ |c| }{$s^{1}_{4}$} &  $a_{1}^{3}$ & $a_{1}^{3}$  & $a_{1}^{1}$ & $a_{1}^{1}$ &
 $a_{1}^{3}$ & $a_{1}^{3}$ & $a_{1}^{1}$ & $a_{1}^{1}$ &
 $a_{1}^{1}$ & $a_{1}^{1}$ & $a_{1}^{1}$ & $a_{1}^{1}$ &
$a_{1}^{1}$ & $a_{1}^{1}$ & $a_{1}^{1}$ & $a_{1}^{1}$
\\
\multicolumn{1}{ |c  }{}                        &
\multicolumn{1}{ |c| }{$s^{1}_{5}$} &  $a_{1}^{3}$ & $a_{1}^{3}$  & $a_{1}^{1}$ & $a_{3}^{32}$ &
 $a_{5}^{12}$ & $a_{5}^{12}$ & $a_{5}^{12}$ & $a_{3}^{23}$ &
 $a_{1}^{1}$ & $a_{4}^{234}$ & $a_{1}^{1}$ & $a_{3}^{23}$ &
$a_{5}^{12}$ & $a_{4}^{234}$ & $a_{5}^{12}$ & $a_{3}^{23}$
\\
\multicolumn{1}{ |c  }{}                        &
\multicolumn{1}{ |c| }{$s^{1}_{6}$} &  $a_{1}^{3}$ & $a_{1}^{3}$  & $a_{1}^{1}$ & $a_{1}^{1}$ &
 $a_{1}^{3}$ & $a_{1}^{3}$ & $a_{1}^{1}$ & $a_{1}^{1}$ &
 $a_{1}^{1}$ & $a_{1}^{1}$ & $a_{1}^{1}$ & $a_{1}^{1}$ &
$a_{1}^{1}$ & $a_{1}^{1}$ & $a_{1}^{1}$ & $a_{1}^{1}$

\\ \cline{1-18}
\multicolumn{1}{ |c  }{\multirow{6}{*}{$s^{3}_{3}$} } &
\multicolumn{1}{ |c| }{$s^{1}_{1}$} &  $a_{2}^{23}$ & $a_{2}^{23}$  & $a_{3}^{1}$ & $a_{2}^{2}$ &
 $a_{5}^{12}$ & $a_{2}^{23}$ & $a_{3}^{1}$ & $a_{2}^{2}$ &
 $a_{4}^{124}$ & $a_{2}^{2}$ & $a_{3}^{1}$ & $a_{2}^{2}$ &
$a_{4}^{124}$ & $a_{2}^{2}$ & $a_{3}^{1}$ & $a_{2}^{2}$
\\
\multicolumn{1}{ |c  }{}                        &
\multicolumn{1}{ |c| }{$s^{1}_{2}$} &  $a_{2}^{3}$ & $a_{2}^{3}$  & $a_{2}^{3}$ & $a_{2}^{3}$ &
 $a_{2}^{3}$ & $a_{2}^{3}$ & $a_{2}^{3}$ & $a_{2}^{3}$ &
 $a_{2}^{3}$ & $a_{2}^{3}$ & $a_{2}^{3}$ & $a_{2}^{3}$ &
$a_{2}^{3}$ & $a_{2}^{3}$ & $a_{2}^{3}$ & $a_{2}^{3}$
\\
\multicolumn{1}{ |c  }{}                        &
\multicolumn{1}{ |c| }{$s^{1}_{3}$} &  $a_{2}^{23}$ & $a_{2}^{23}$  & $a_{23}^{2}$ & $a_{2}^{2}$ &
 $a_{5}^{1}$ & $a_{2}^{23}$ & $a_{5}^{1}$ & $a_{2}^{2}$ &
 $a_{4}^{14}$ & $a_{2}^{2}$ & $a_{4}^{14}$ & $a_{2}^{2}$ &
$a_{4}^{14}$ & $a_{2}^{2}$ & $a_{4}^{14}$ & $a_{2}^{2}$
\\
\multicolumn{1}{ |c  }{}                        &
\multicolumn{1}{ |c| }{$s^{1}_{4}$} &  $a_{2}^{3}$ & $a_{2}^{3}$  & $a_{2}^{3}$ & $a_{2}^{3}$ &
 $a_{2}^{3}$ & $a_{2}^{3}$ & $a_{2}^{3}$ & $a_{2}^{3}$ &
 $a_{2}^{3}$ & $a_{2}^{3}$ & $a_{2}^{3}$ & $a_{2}^{3}$ &
$a_{2}^{3}$ & $a_{2}^{3}$ & $a_{2}^{3}$ & $a_{2}^{3}$
\\
\multicolumn{1}{ |c  }{}                        &
\multicolumn{1}{ |c| }{$s^{1}_{5}$} &  $a_{2}^{23}$ & $a_{2}^{23}$  & $a_{2}^{23}$ & $a_{2}^{2}$ &
 $a_{5}^{1}$ & $a_{2}^{23}$ & $a_{5}^{1}$ & $a_{2}^{2}$ &
 $a_{2}^{23}$ & $a_{2}^{2}$ & $a_{2}^{23}$ & $a_{2}^{2}$ &
$a_{5}^{1}$ & $a_{2}^{2}$ & $a_{5}^{1}$ & $a_{2}^{2}$
\\
\multicolumn{1}{ |c  }{}                        &
\multicolumn{1}{ |c| }{$s^{1}_{6}$} &  $a_{2}^{3}$ & $a_{2}^{3}$  & $a_{2}^{3}$ & $a_{2}^{3}$ &
 $a_{2}^{3}$ & $a_{2}^{3}$ & $a_{2}^{3}$ & $a_{2}^{3}$ &
 $a_{2}^{3}$ & $a_{2}^{3}$ & $a_{2}^{3}$ & $a_{2}^{3}$ &
$a_{2}^{3}$ & $a_{2}^{3}$ & $a_{2}^{3}$ & $a_{2}^{3}$
\\ \cline{1-18}
\multicolumn{1}{ |c  }{\multirow{6}{*}{$s^{3}_{4}$} } &
\multicolumn{1}{ |c| }{$s^{1}_{1}$} &  $a_{1}^{3}$ & $a_{2}^{23}$  & $a_{3}^{2}$ & $a_{2}^{2}$ &
 $a_{5}^{12}$ & $a_{2}^{23}$ & $a_{3}^{2}$ & $a_{2}^{2}$ &
 $a_{4}^{24}$ & $a_{2}^{2}$ & $a_{3}^{1}$ & $a_{2}^{2}$ &
$a_{4}^{24}$ & $a_{2}^{2}$ & $a_{3}^{1}$ & $a_{2}^{2}$
\\
\multicolumn{1}{ |c  }{}                        &
\multicolumn{1}{ |c| }{$s^{1}_{2}$} &  $a_{1}^{3}$ & $a_{1}^{13}$  & $a_{1}^{1}$ & $a_{1}^{1}$ &
 $a_{1}^{3}$ & $a_{1}^{13}$ & $a_{1}^{1}$ & $a_{1}^{1}$ &
 $a_{1}^{1}$ & $a_{1}^{1}$ & $a_{1}^{1}$ & $a_{1}^{1}$ &
$a_{1}^{1}$ & $a_{1}^{1}$ & $a_{1}^{1}$ & $a_{1}^{1}$
\\
\multicolumn{1}{ |c  }{}                        &
\multicolumn{1}{ |c| }{$s^{1}_{3}$} &  $a_{1}^{3}$ & $a_{2}^{23}$  & $a_{1}^{1}$ & $a_{2}^{2}$ &
 $a_{5}^{12}$ & $a_{2}^{23}$ & $a_{5}^{12}$ & $a_{2}^{2}$ &
 $a_{4}^{4}$ & $a_{2}^{2}$ & $a_{4}^{4}$ & $a_{2}^{2}$ &
$a_{4}^{4}$ & $a_{2}^{2}$ & $a_{4}^{4}$ & $a_{2}^{2}$
\\
\multicolumn{1}{ |c  }{}                        &
\multicolumn{1}{ |c| }{$s^{1}_{4}$} &  $a_{1}^{3}$ & $a_{1}^{13}$  & $a_{1}^{1}$ & $a_{1}^{1}$ &
 $a_{1}^{3}$ & $a_{1}^{13}$ & $a_{1}^{1}$ & $a_{1}^{1}$ &
 $a_{1}^{1}$ & $a_{1}^{1}$ & $a_{1}^{1}$ & $a_{1}^{1}$ &
$a_{1}^{1}$ & $a_{1}^{1}$ & $a_{1}^{1}$ & $a_{1}^{1}$
\\
\multicolumn{1}{ |c  }{}                        &
\multicolumn{1}{ |c| }{$s^{1}_{5}$} &  $a_{1}^{3}$ & $a_{2}^{23}$  & $a_{1}^{1}$ & $a_{2}^{2}$ &
 $a_{5}^{12}$ & $a_{2}^{23}$ & $a_{5}^{12}$ & $a_{2}^{2}$ &
 $a_{1}^{1}$ & $a_{2}^{2}$  & $a_{1}^{1}$ & $a_{2}^{2}$ &
 $a_{5}^{12}$ & $a_{2}^{2}$ & $a_{5}^{12}$ & $a_{2}^{2}$
\\
\multicolumn{1}{ |c  }{}                        &
\multicolumn{1}{ |c| }{$s^{1}_{6}$} &  $a_{1}^{3}$ & $a_{1}^{13}$  & $a_{1}^{1}$ & $a_{1}^{1}$ &
 $a_{1}^{3}$ & $a_{1}^{13}$ & $a_{1}^{1}$ & $a_{1}^{1}$ &
 $a_{1}^{1}$ & $a_{1}^{1}$ & $a_{1}^{1}$ & $a_{1}^{1}$ &
$a_{1}^{1}$ & $a_{1}^{1}$ & $a_{1}^{1}$ & $a_{1}^{1}$
\\ \cline{1-18}
\end{tabular}

\caption{
This table represents our main example in the normal form.
\newline
The game form  $g : S \rightarrow A$, in which
$S= S_1 \times S_2 \times S_3 \times S_4$  and  $A = \{a_1, a_2, a_3, a_4, a_5; c\}$,
is given by the four-dimensional table of size  $6 \times 8 \times 4 \times 2$.
\newline
Player $1$ has six strategies:
%%% $s^{1}_{i} (u_1, v_1)$ are:
\newline
$s^{1}_{1}$: $(u_1, u_2), (v_1, v_2)$,
$s^{1}_{2}$: $(u_1, u_3), (v_1, v_2)$,
\newline
$s^{1}_{3}$: $(u_1, u_2), (v_1, u_4)$,
$s^{1}_{4}$: $(u_1, u_3), (v_1, u_4)$,
\newline
$s^{1}_{5}$: $(u_1, u_2), (v_1, w_2)$,
$s^{1}_{6}$: $(u_1, u_3), (v_1, w_2)$;
\newline
player $2$ has eight strategies:
\newline
%%% $s^{2}_{i}(u_2, v_2, w_2)$ are:
$s^{2}_{1}$: $(u_2, v_1), (v_2, u_4), (w_2, u_3)$,
$s^{2}_{2}$: $(u_2, v_3), (v_2, u_4), (w_2, u_3)$,
\newline
$s^{2}_{3}$: $(u_2, v_1), (v_2, a_3), (w_2, u_3)$,
$s^{2}_{4}$: $(u_2, v_3), (v_2, a_3), (w_2, u_3)$,
\newline
$s^{2}_{5}$: $(u_2, v_1), (v_2, u_4), (w_2, a_5)$,
$s^{2}_{6}$: $(u_2, v_3), (v_2, u_4), (w_2, a_5)$,
\newline
$s^{2}_{7}$: $(u_2, v_1), (v_2, a_3), (w_2, a_5)$,
$s^{2}_{8}$: $(u_2, v_3), (v_2, a_3), (w_2, a_5)$;
\newline
player $3$ has four strategies:
%%% $s^{3}_{i}(u_3, v_3)$ are:
\newline
$s^{3}_{1}$: $(u_3, v_3), (v_3, v_2)$,
$s^{3}_{2}$: $(u_3, a_1), (v_3, v_2)$,
$s^{3}_{3}$: $(u_3, v_3), (v_3, a_2)$,
$s^{3}_{4}$: $(u_3, a_1), (v_3, a_2)$;
\newline
finally, player $4$ has two strategies:
$s^{4}_{1}(u_4, w_2)$,
$s^{4}_{2}(u_4, a_4)$.
\newline
It is not difficult
(although time consuming) to verify that
the corresponding game has no NE
whenever the preference profile
$o = \{o_1, o_2, o_3, o_4\}$  of the players
$I = \{1,2,3,4\}$  agrees with the following partial orders:
\newline
$O_1: a_2 > a_4 > a_3 > a_1 > a_5; \;\;$
% \newline
$O_2 : \min(a_1, c) > a_3 > \max(a_4, a_5) > \min(a_4, a_5) > a_2; \;\;$
\newline
$O_3 : \min(a_5, c) > a_1 > a_2 > \max(a_3, a_4); \;\;$
%\newline
$O_4 : \min(a_1, a_2, a_3, a_5) > a_4 > c$.
\newline
For every situation
$s = (s^1_{\ell_1}, s^2_{\ell_2}, s^3_{\ell_3}, s^4_{\ell_4}) \in S_1 \times S_2 \times S_3 \times S_4 = S$
the outcome  $g(s)$, which is either a terminal  $a_j$  or the dicycle  $c$,  is
shown in the entry  $(\ell_1, \ell_2, \ell_3, \ell_4)$  of the table.
The upper indices indicate the players who can improve the situation  $s$.
Thus, a situation  $s$  is a NE if and only if the corresponding outcome has no upper indices.
Since the table contains no such situation, the considered game has no NE.
}
\end{table}

\section{The main example}
\label{s1}

The positional and normal forms
of the game announced in the title of the paper are
presented below by Figure \ref{f1} and Table \ref{t1}, respectively.

\newpage

\section{Open ends}
\label{s2}

In the above example there are four players, $n=4$,
five terminals, $p=5$, and there is only one dicycle.

\smallskip

In \cite{BG03} it was shown that
there is no chess-like NE-free game with  $p \leq 2$  terminals.
In \cite{BR09}, Boros and Rand extended this result $p \leq 3$.
Thus, only the case  $p=4$ remains open.

\smallskip

It is known that
every two-person chess-like game has a NE;
see \cite{BG03, BGMS07}, the last section of each paper.
Thus, only the case  $n=3$  remains open.

\smallskip

The proof for  $n=2$  is simple and we
repeat it here for convenience of the readers.
It is based on the following important property
of the two-person game forms, which
seems not to be extendable for  $n > 2$.

A two-person game form  $g$  is called
(i) {\em Nash-solvable},
(ii) {\em zero-sum-solvable}, and
(iii) {\em $\pm 1$-solvable}  if
the corresponding game  $(g,u)$  has at least one NE
(i) for every payoff   $u = (u_1, u_2); \;\;$
(ii) for every payoff   $u = (u_1, u_2)$  such that
$u_1(a) + u_2(a) = 0$  for each outcome  $a \in A; \;\;$
(iii) for every payoff  $u = (u_1, u_2)$  such that
$u_1(a) + u_2(a) = 0$  for each outcome  $a \in A$
and both  $u_1$  and  $u_2$  take only values  $+1$  or  $-1$.

In fact, all three above properties of a game forms are equivalent.
For (ii) and (iii) this was shown in 1970 by
Edmonds and Fulkerson \cite{EF70} and independently in \cite{Gur73}.
Then, the list was extended by statement (i)
in \cite{Gur75}; see also \cite{Gur88}, where
it also shown that
a similar statement fails for the three-person game forms.

Thus, it is sufficient to prove $\pm 1$-solvability,
rather than  Nash-solvability, of the two-person chess-like games.
Hence, we can assume that
each outcome  $a \in A = V_T \cup \{c\}$
is either winning for player
$1$  and losing for player  $2$,  or vice versa.
Without any loss of generality, assume
that  $c$  is winning for  $1$.

Then, let  $V_T = V_T^1 \cup V_T^2$
be the partition of all terminals into
outcomes winning for players  $1$ and $2$, respectively.
Furthermore, let  $V^2 \subseteq V$  denote the set of all positions
from which player  $2$  can enforce  $V_T^2$;
in particular, $V_T^2 \subseteq V^2$.
Finally, let us set   $V^1 = V \setminus V^2$;
in particular, $V_T^1 \subseteq V^1$.

By the above definitions, in every position
$v \in V_1 \cap V^1$  player  $1$  can stay
out of   $V^2$, that is, (s)he has a move
$(v, v')$  such that  $v' \in V^1$.
Let us fix a strategy  $s^1_0$  that chooses such a move
in each position  $v \in V_1 \cap V^1$
and any move in  $v \in V_1 \cap V^2$.
Then,  for any  $s^2 \in S_2$, the outcome  $g(s^1_0, s^2)$
is winning for player  $1$  whenever
the initial position  $v_0$  is in  $V^1$.
Indeed, either  $g(s^1_0, s_2) \in V^1_T$, or $g(s^1_0, s^2) = c$;
in both cases player  $1$  wins.
Thus, player $1$ wins when  $v_0 \in V^1$ and
player $2$  wins when  $v_0 \in V^2$;
in each case a saddle point exists.

\medskip

In \cite{GO14},  we conjectured that
there is a NE-free chess-like game
satisfying the following extra condition

\begin{itemize}
\item{(C)} outcome  $c$  is worse than
any terminal outcome $a \in V_T$ for each player $i \in I$.
\end{itemize}

This conjecture remains open.
% It is important, because
Such an example, if exists, would
strengthen simultaneously the example of the present paper and
the main example of \cite{GO14}; see Figure 2 and Table 2 there.

\begin{remark}
In \cite{BEGM11}, $16$ problems related to Nash-solvability, subgame perfect and
not, of the chess-like and backgammon-like games were
considered, under the assumption  (C)  and without it.
For $15$ of these problems, condition (C)  appears ``irrelevant", that is,
either Nash-solvability holds, even without  (C), or
even with  (C), it fails.  %%%  hold,  even if  (C)  is assumed.
(Yet, in the latter case, the size of the corresponding example may increase significantly.)
Based on these observations, we conjectured in \cite{GO14} that
the same will happen for the last $16$th case, which is the subject of the present paper.
An example might be similar to one in Figure \ref{f1}, but with a much larger digraph.
Several interpretations of assumption  (C)  were suggested
Remark 2 and Section 2.1.1 of  \cite{BG03}.
\end{remark}

Finally, it follows from the main result of \cite{BG03} that
there is no NE-free  $n$-person chess-like game in which
each player controls a unique position.
In the above example, the players  $1,2,3,$ and $4$  control
$2,3,2,$ and $1$  positions, respectively.
It remains open, if there is a chess-like NE-free game
in which each player controls, say, at most two positions.

\section{Related results on Nash-solvability}
\label{s3}

In the next two subsections we recall
two large families of games with perfect information that
are known to be Nash-solvable in pure stationary uniformly optimal strategies.

\subsection{Acyclyc $n$-person
backgammon-like games with perfect information}
\label{ss20}

In 1950 Nash introduced his concept of equilibrium
for the normal form  $n$-person games \cite{Nas50, Nas51}.
Soon after, Kuhn \cite{Kuh50, Kuh53} and Gale \cite{Gal53}
suggested the backward induction procedure and proved that
any finite acyclic chess-like game with perfect information has
a NE in pure stationary strategies;
moreover, the obtained NE is {\em subgame perfect}, that is, the same
strategy profile is a NE with respect to any initial position.
The authors restricted themselves
to the chess-like games on finite arborescence
(directed trees) but, in fact, backward induction can be easily
extended to the backgammon-like games on the finite digraphs
without dicycles (so called {\em acyclic digraphs or DAGs}).
Yet, acyclicity is a crucial assumption and cannot be waved.

\smallskip

For any integer  $k \geq 2$  let us introduce a digraph  $G_k$  that
consists of  $k$  terminals  $a_1, \ldots, a_k$,
the directed  $k$-cycle
$C_k$  on the  $k$  non-terminal vertices  $v_1, \ldots, v_k$, and
the perfect matching  $(v_j, a_j), j = 1, \ldots, k$
between these vertices and the terminals.

The existence of a subgame perfect NE fails already for  $k=2$  \cite{AGH10}.
Let players  $1$ and $2$ control vertices  $v_1$  and $v_2$  and
have the preferences, $o_1: (c > a_1 > a_2)$  and  $(a_1 > a_2 > c)$, respectively.
It is easy to verify that a NE exists for any given initial position, $v_1$  or  $v_2$, but
no strategy profile is a NE with respect to both simultaneously.
Let us notice that
the preferences are not opposite
(both players prefer  $a_1$  to  $a_2$), while
$c$  is the worst outcome for player  $2$  and the best one for  $1$.

A similar example exists even if in addition we require  (C):
the dicycle is worse than each terminal for both players.
Consider digraph  $G_6$  in which players  $1$ ad $2$ control
the odd and even positions
($a_1, a_3, a_5$  and  $a_2, a_4, a_6$) respectively.
It was shown in \cite{BEGM11} that there exists no subgame perfect NE whenever

\medskip
\noindent
$o_1: a_6 > a_5 > a_2 > a_1 > a_3 > a_4 > c;$  % and
$o_2 \in O_2: \{a_3 > a_2 > a_6 > a_4 > a_5 > c, a_6 > a_1 > c\}$.

\medskip

Note that there is no such example for $G_4$.
A similar example for a three-person game was constructed in \cite{BEGM11}.
The players $1,2,3$  control, respectively,
the positions $v_1, v_2, v_3$  of  $G_3$  and have the preferences:

\medskip
\noindent
$o_1: a_2 > a_1 > a_3 > c; \;\;\;$
$o_2: a_3 > a_2 > a_1 > c; \;\;\;$
$o_3: a_1 > a_3 > a_2 > c.$

\smallskip

In other words, for all players: $c$  is the worst outcome, in accordance with (C);
it is better if the previous player terminate; it is still better to terminate himself;
finally, the best if the next player terminates.

\medskip

In \cite{BGY13}, these results were strengthen as follows.
It was demonstrated that a subgame perfect NE may fail to exist
not only in the pure but even in thee mixed strategies.
The corresponding examples are based on the same
positional game forms, $G_6$ for $n=2$ and $G_3$ for $n=3$, but
the above preference profiles
% $o_1, o_2$ and respectively  $o_1, o_2, o_3$
are replaced by some carefully chosen payoffs, which agree
with the corresponding preferences.

\bigskip

The above examples imply that, for  any  $n \geq 2$,  an $n$-person
{\em backgammon-like} game, even with a fixed initial position, may have no NE.
Given a chess-like game  $(G = (V, E), D, o)$, in which
no initial position is fixed yet,
add to it a new position  $v_0$  and
the move  $(v_0, v)$  from  $v_0$
to each non-terminal position $v \in V \setminus V_T$.
Furthermore, assign to  $(v_0, v)$
a non-negative probability  $p_v \geq 0$  such that  $\sum_{v \in V}  p_v = 1$.
Denote by  $(G', D, o, v_0)$  the obtained
backgammon-like game form, in which  $v_0$  is the initial position;
see Figure 1 in \cite{BEGM11}. It was shown in \cite{BEGM11} that

\begin{itemize}
%\noindent
\item{(i)} if  $s$  is a subgame perfect NE in  $(G, D, o)$
then  $s$  is a NE in  $(G', D, o, v_0)$  for any  $p_v$;
%\noindent
\item{(ii)} if  $s$  is a  NE in  $(G', D, o, v_0)$
and  $p_v > 0  \, \forall v \in V \setminus V_T$, then
$s$  is a subgame perfect NE in  $(G_k, D, o)$.
\end{itemize}

These results imply that
Nash-solvability of  $(G', D, v_0)$  is equivalent with
subgame perfect Nash-solvability of  $(G, D)$.
As we know, the latter property may fail for  $G = G_k$  for any  $n \geq 2$.
Thus,  for any  $n \geq 2$, an  $n$-person backgammon-like game
$(G'_k, D, v_0, o)$,  in which $v_0$  is a unique position of chance, may have no NE.

\subsection{Two-person zero-sum % chess- and backgammon-like
games with perfect information}
\label{ss21}

According to the previous subsection, in
the presence of dicycles, backward induction fails, in general.
Yet, it can be modified
(and thus saved) in case of the two-person zero-sum games
(that is, when  $I = \{1,2\}$  and
$u(1, a) + u(2, a) = 0$  for any outcome $a \in A$,
or two preferences  $o_1$  and  $o_2$
of the players  $1$ and $2$  over  $A$  are opposite).

For example, the recent paper \cite{GGABR14} shows
how to solve by backward induction
a two-person zero-sum game that is ``acyclic", but the players can pass;
in other words, the corresponding digraph
contains a loop at each vertex, but no other dicycles.

A general linear time algorithm solving
ay two-person zero-sum chess-like game,
by a modified backward induction, was suggested in \cite{AHMS10}
and independently in \cite{BEGM11}.
In contrast, no polynomial algorithm is known for the
two-person zero-sum backgammon-like games \cite{Con92, Con93}.
However, it is well-known that
subgame perfect saddle points in stationary strategies exist
in this case and even in much more general cases considered below.

In fact, studying two-person zero-sum chess-like games began
long before the backward induction was suggested in early fifties by
\cite{Kuh50, Kuh53, Gal53}.
Zermelo gave his seminal talk on solvablity of chess
in pure strategies \cite{Zer1912} as early as in 1912.
Later, K\"onig  \cite{Kon27} and Kalman \cite{Kal28}
strengthen this result showing that
there exist pure stationary uniformly optimal strategies
producing a subgame perfect saddle point, in
any two-person zero-sum chess-like game.
The reader can find  more detailed surveys in
\cite{Was90, Ewe01, SW01, AHMS10};
see also \cite{BG09, FKSV09, FKSV10, Gur90a}.

The chess-like and backgammon-like games, considered
in this note, by the definition, are {\em transition-free}.
There is a much more general class:
{\em stochastic games with perfect information} in which
a {\em transition payoff}  $r_i(u,v)$  is defined
for every move  $(u,v)$  and for each player  $i \in I$.
Gillette, in his seminal paper \cite{Gil57},
introduced the {\em mean (or average)} effective payoff
for these  games and proved the existence of
subgame perfect saddle point in the uniformly optimal stationary strategies
for the two-person zero-sum case.
The proof is pretty complicated.
It is based on the Tauberian theory
and, in particular, on the Hardy-Littlewood theorem \cite{HL31}.
(In \cite{Gil57}, the conditions of this theorem were not accurately
verified and the flaw was corrected in twelve years by
Liggett and Lippman in \cite{LL69}.)

Stochastic games with perfect information can be viewed as
{\em backgammon-like} games with transition payoffs.
More precisely, these two classes are polynomially equivalent \cite{BEGM13}.
Interestingly, the corresponding two-person zero-sum {\em chess-like} games
(with transition payoffs but without random moves) so-called
{\em cyclic mean-payoff games}, appeared only 20-30 years later, introduced 
for the complete bipartite digraphs by Moulin \cite{Mou76, Mou76a},
for any bipartite digraphs by  Ehrenfeucht and Mycielski
\cite{EM73, EM79},  and for arbitrary digraphs by
Gurvich, Karzanov, and Khachiyan \cite{GKK88}.
Again, the existence of a saddle point in
the pure stationary uniformly optimal strategies
was proven for the two-person zero-sum case.

This result cannot be extend to the non-zero-sum case.
In \cite{Gur88}, a cyclic mean payoff two-person NE-free game
was constructed on the complete bipartite  $3 \times 3$  digraph
with symmetric payoffs. 
(The corresponding normal form game is a $27 \times 27$  bimatrix.)
It was shown in \cite{Gur90} that
this example is, in a sense, minimal, namely, a NE always exists
for the games on the complete  $(2 \times \ell)$  bipartite digraphs.

\smallskip

A general family of the so-called  $k$-{\em total effective payoffs} was
recently introduced in \cite{BEGM14-k} for any nonnegative integer $k$ such that
the $0$-total one is the mean effective payoff, while
the $1$-total one is the total effective payoff
introduced earlier by Thuijsman and Vrieze in \cite{TV87, TV98}.
The existence of a saddle point in uniformly optimal pure stationary strategies
for the two-person zero-sum chess-like games 
with the $k$-total effective payoff was proven for all  $k$  in \cite{BEGM14-k}.
If this result can be extended to the backgammon-like games is an open problem.
Yet, for  $k \leq 1$  the answer is positive.
As was already mentioned, for $k=0$  it was proven long ago.
For $k=1$  the result was first obtained in \cite{TV87, TV98}, see also \cite{BEGM14-k}.

However, it cannot be extended to the non-zero-sum case.
In particular, a NE (in pure stationary strategies) may fail to exist 
already in a two-person chess-like game.
For  $k=0$  the example was given in \cite{Gur88}.
Furthermore, in \cite{BEGM14-k}, a simple embedding
of the  $(k-1)$-total payoff games into
the $k$-total ones was constructed.
Thus, the example of  \cite{Gur88}  works for all  $k$.

\bigskip

{\bf Acknowledgemets}: The authors are thankful to
Endre Boros, Konrad Borys, Khaled Elbassioni, and Gabor Rudolf
for helpful discussions.

\end{document}